# Global uniform asymptotic stabilization and $\kappa-$ exponential trajectory tracking of underactuated surface ships with non-diagonal inertia/damping matrices

Ma Bao-Li

The Seventh Research Division, BeiHang University, Beijing 100191

Email:mabaoli@buaa.edu

**Abstract:** In this work, we investigate the state stabilization and trajectory tracking problems of underactuated surface ships with full state model of having non-diagonal inertia and damping matrices. By combining the novel state transformations, the direct Lyapunov approach, and the nonlinear time-varying tools, the stabilization and the trajectory tracking controllers are developed respectively guaranteeing global uniform asymptotic convergence of the state to the desired set point and global $\kappa-$ exponential convergence to the desired reference trajectory via mild persistent exciting conditions. Simulation examples are given to illustrate the effectiveness of the proposed control schemes.

## 1. Introduction

Control of surface ships has been an active research topic in recent years [1]. With the simplifying full ship model of having diagonal inertia and damping matrices, various approaches have been developed for point stabilization [2-8], trajectory tracking [8-12], both the point stabilization and trajectory tracking [13-14], and path following [15-17]. Realizing that such a simplifying model is unrealistic, recent research has aimed at dealing with the full ship model without the simplifying assumption of having diagonal inertia and damping matrices [18-22]. However, the presented results are limited to the problems of path following [18-21] and position control [22]. As our best knowledge, there exists no control scheme dealing with full state stabilization and full state trajectory tracking for the ship model with non-diagonal inertia/damping matrices even in the absence of model uncertainties and external disturbances.

The purpose of this work is to solve the point stabilization and trajectory tracking problems of surface ships with non-zero off-diagonal terms presented in their inertia and damping matrices. By introducing novel coordinate and input transformations, the error models for stabilization and trajectory tracking are converted to the new cascaded forms respectively, which allow for deriving a smooth time-varying control law able to achieve global uniform asymptotic full state point stabilization, and allow for obtaining a continuous control law able to achieve global $\kappa-$exponential trajectory tracking via mild persistent exciting (PE) conditions imposed on the reference trajectory.

The reminder of the work is organized as follows. After introducing the full state model of underactuated surface ships with non-diagonal inertia/damping matrices, the point stabilization and trajectory tracking problems are formulated in Section 2. The point stabilization controller is developed in Section 3, and the trajectory tracking one is proposed in Section 4. Simulation examples are illustrated in Section 5. Section 6 concludes the work.

**Notations**：

A continuous function $\rho(r):[0,a) \to [0,\infty)$ is said to belong to class-$\kappa$ if it is strictly increasing and $\rho(0)=0$. It is said to belong to class $\kappa_\infty$ if $a=\infty$ and $\rho(r) \to \infty$ as $r \to \infty$ [23].

A continuous function $\chi(r,s):[0,a)\times[0,\infty) \to [0,\infty)$ is said to belong to class $\kappa\ell$ if, for each fixed $s$, the mapping $\chi(r,s)$ belongs to class $\kappa$ with respect to $r$ and, for each fixed $r$, the mapping $\chi(r,s)$ is decreasing with respect to $s$ and $\chi(r,s) \to 0$ as $s \to \infty$ [23].

A system $\dot{x}=f(t,x), f(t,0)=0, \forall t \geq t_0$ (with $x \in R^n$ and $f(t,x)$ piecewise continuous in $t$ and locally Lipschitz in $x$) is called global uniform asymptotic stable (**GUAS**) if there exists a class $\kappa\ell$ function $\chi(\cdot,\cdot)$ such that the state trajectory starting from any initial state $x(t_0) \in R^n$ satisfy $\|x(t)\| \leq \chi(\|x(t_0)\|, t-t_0) e^{-\gamma(t-t_0)}$ [23]. It is called global $\kappa$-exponential stable (**GKES**) if there exist a class $\kappa$ function $\rho(\cdot)$ and a positive constant $\gamma$ such that the state trajectory starting from any initial state $x(t_0) \in R^n$ satisfy $\|x(t)\| \leq \rho(\|x(t_0)\|) e^{-\gamma(t-t_0)}$ [24].

2. **Problem formulation**

We consider the ship model described by [1]

$$\mathbf{M}\dot{\mathbf{v}} + \mathbf{C}(\mathbf{v})\mathbf{v} + \mathbf{D}\mathbf{v} = \boldsymbol{\tau} \quad (1)$$

where $\mathbf{v} = [u, v, r]^T$ denotes the velocities in surge, sway and yaw directions respectively, $\boldsymbol{\tau} = [\tau_u, 0, \tau_r]^T$ denotes the control input with $(\tau_u, \tau_r)$ the control force and moment in surge and yaw directions respectively. The system matrices are given by

$$\mathbf{M} = \begin{bmatrix} m_{11} & 0 & 0 \\ 0 & m_{22} & m_{23} \\ 0 & m_{23} & m_{33} \end{bmatrix},$$

$$\mathbf{C}(\mathbf{v}) = \begin{bmatrix} 0 & 0 & -m_{22}v - m_{23}r \\ 0 & 0 & m_{11}u \\ m_{22}v + m_{23}r & -m_{11}u & 0 \end{bmatrix}, \quad (2)$$

$$\mathbf{D} = \begin{bmatrix} d_{11} & 0 & 0 \\ 0 & d_{22} & d_{23} \\ 0 & d_{32} & d_{33} \end{bmatrix}$$

where $d_{11} > 0, d_{22} > 0, d_{23} \neq 0, d_{33} > 0, m_{11} > 0, m_{22} > 0, m_{23} \neq 0, m_{33} > 0$ represent the hydrodynamic damping and ship's inertia including added mass.

The kinematics of the ship is described by

$$\begin{bmatrix} \dot{x} \\ \dot{y} \\ \dot{\psi} \end{bmatrix} = \begin{bmatrix} \cos\psi & -\sin\psi & 0 \\ \sin\psi & \cos\psi & 0 \\ 0 & 0 & 1 \end{bmatrix} \begin{bmatrix} u \\ v \\ r \end{bmatrix} \quad (3)$$

where $x$, $y$ and $\psi$ represent the position and orientation of the ship in the Earth-fixed frame.

In modeling the vessel dynamics, it is assumed that: A1) the high-order damping terms and the dynamics associated with the motions in heave, roll and pitch are all negligible; A2) the inertia and damping matrices $(\mathbf{M}, \mathbf{D})$ are known; A3) the external disturbances are not included. The stabilization and tracking problems by relaxing these assumptions will be investigated in the future research.

Based on the standard feedback linearization procedure used in [20], the dynamic equation (1) can be simplified to

$$\begin{cases} \dot{u} = \tau_1, \\ \dot{r} = \tau_2, \\ \dot{v} = \dfrac{1}{m_{22}}(-m_{23}\dot{r} - m_{11}ur - d_{22}v - d_{23}r) \\ \quad = -a\tau_2 - br - cur - dv \end{cases} \quad (4)$$

where $a = \dfrac{m_{23}}{m_{22}} \neq 0, b = \dfrac{d_{23}}{m_{22}} \neq 0, c = \dfrac{m_{11}}{m_{22}} > 0, d = \dfrac{d_{22}}{m_{22}} > 0$, and $(\tau_1, \tau_2)$ are the new control inputs, which are related to the true ones $(\tau_u, \tau_r)$ via the following input transformation

$$\begin{aligned} \tau_1 &= \dfrac{1}{m_{11}}\left(\tau_u - r(-m_{22}v - m_{23}r) - d_{11}u\right), \\ \tau_2 &= \dfrac{1}{m_{22}m_{33} - m_{23}^2}\begin{bmatrix} -m_{23} & m_{22} \end{bmatrix}\begin{bmatrix} -(m_{11}ur + d_{22}v + d_{23}r) \\ \tau_r - ((m_{22}v + m_{23}r)u - m_{11}uv + d_{32}v + d_{33}r) \end{bmatrix} \\ &= \dfrac{1}{m_{22}m_{33} - m_{23}^2}\left(m_{22}\tau_r + (m_{11} - m_{22})(m_{23}r + m_{22}v)u + (m_{23}d_{22} - m_{22}d_{32})v + (m_{23}d_{23} - m_{22}d_{33})r\right) \end{aligned} \quad (5)$$

The input transformation (5) is clearly global invertible provided $m_{22}m_{33} - m_{23}^2 > 0$, which is satisfied for any positive symmetric definite matrix $\mathbf{M}$.

The inverse of the input transformation (5) can be derived as

$$\begin{aligned} \tau_u &= m_{11}\tau_1 - r(m_{22}v + m_{23}r) + d_{11}u, \\ \tau_r &= \dfrac{(m_{22}m_{33} - m_{23}^2)\tau_2 - \left((m_{11} - m_{22})(m_{23}r + m_{22}v)u + (m_{23}d_{22} - m_{22}d_{32})v + (m_{23}d_{23} - m_{22}d_{33})r\right)}{m_{22}} \end{aligned} \quad (6)$$

The control object for point stabilization of surface ships can be stated as: *find a feedback control law* $\tau_1(\cdot), \tau_2(\cdot)$ *such that the origin of the closed-loop system (3)-(4) is GUAS.*

For trajectory tracking, the reference trajectory is generated by the same model of (3)-(4) as follows

$$\begin{bmatrix} \dot{x}_d \\ \dot{y}_d \\ \dot{\psi}_d \end{bmatrix} = \begin{bmatrix} \cos\psi_d & -\sin\psi_d & 0 \\ \sin\psi_d & \cos\psi_d & 0 \\ 0 & 0 & 1 \end{bmatrix}\begin{bmatrix} u_d \\ v_d \\ r_d \end{bmatrix} \quad (7)$$

$$\begin{cases} \dot{u}_d = \tau_{1d}, \\ \dot{r}_d = \tau_{2d}, \\ \dot{v}_d = -a\tau_{2d} - br_d - cu_d r_d - dv_d \end{cases} \quad (8)$$

The control object for global $\kappa$ exponential trajectory tracking can be stated as: *find a feedback control law* $\tau_1(\cdot), \tau_2(\cdot)$ *such that the state tracking errors* $(x - x_d, y - y_d, \psi - \psi_d,$

$u - u_d, v - v_d, r - r_d$) *globally asymptotically converge to zero with an exponential convergence rate.*

## 3. Point stabilization

In this section, the point stabilization scheme is derived by introducing novel coordinate and input transformations, and the stability of the closed-loop system is analyzed by Lyapunov approach and nonlinear time-varying tools.

The first coordinate transformation is the same as that introduced in [2]

$$\bar{x} = x\cos\psi + y\sin\psi, \bar{y} = -x\sin\psi + y\cos\psi \qquad (9)$$

The dynamics of $(\bar{x}, \bar{y})$ coordinates is obtained as

$$\dot{\bar{x}} = u + r\bar{y}, \dot{\bar{y}} = v - r\bar{x}$$

The second coordinate and input transformations are as follows

$$\begin{aligned}
\bar{v} &= v + ar + b\psi, \\
z &= d\bar{y} + \bar{v}, \\
\bar{u} &= cu + d\bar{x}, \\
\bar{\tau}_1 &= c\tau_1 + d\frac{\bar{u} - d\bar{x}}{c} + (z - \bar{v})r
\end{aligned} \qquad (10)$$

which is novel and introduced to remove the term $-a\tau_2 - br$ in the sway dynamics (third equation of (4)) and the term of $v$ in the $\bar{y}$ dynamics such that the ship's model can be converted to the following advantageous form

$$\begin{aligned}
\dot{\bar{x}} &= u + r\bar{y} = \frac{\bar{u} - d\bar{x}}{c} + \frac{z - \bar{v}}{d}r, \\
\dot{\bar{v}} &= -a\tau_2 - br - cur - dv + a\tau_2 + br = -(\bar{u} - d\bar{x})r - d(\bar{v} - ar - b\psi), \\
\dot{z} &= dv - dr\bar{x} - (\bar{u} - d\bar{x})r - d(\bar{v} - ar - b\psi) = -\bar{u}r, \\
\dot{\psi} &= r, \\
\dot{\bar{u}} &= c\tau_1 + d\dot{\bar{x}} = c\tau_1 + d\frac{\bar{u} - d\bar{x}}{c} + (z - \bar{v})r \triangleq \bar{\tau}_1, \\
\dot{r} &= \tau_2
\end{aligned} \qquad (11)$$

Equation (11) can be rewritten into the following two subsystems

$$\dot{\bar{x}} = -\frac{d}{c}\bar{x} - \frac{1}{d}\bar{v}r + D_1, \qquad (12)$$
$$\dot{\bar{v}} = -d\bar{v} + d\bar{x}r + D_2$$

$$\dot{z} = -\bar{u}r,$$
$$\dot{\psi} = r,$$
$$\dot{\bar{u}} = \bar{\tau}_1, \qquad (13)$$
$$\dot{r} = \tau_2$$

where

$$D_1 = \frac{\bar{u}}{c} + \frac{zr}{d}, \quad D_2 = -\bar{u}r + d(ar + b\psi) \qquad (14)$$

**Remark 1.** The novel coordinate and input transformations (10) are crucial to convert the ship's model into the cascaded form (12)-(13) such that the stabilization problem of the full ship model (12)-(13) is reduced to the one of its subsystem (13) as claimed in the following Lemma.

**Lemma 1.** Any control law $\bar{\tau}_1(t,z,\psi,\bar{u},r), \tau_2(t,z,\psi,\bar{u},r)$ that makes the subsystem (13) GUAS also makes the whole system (12)-(13) GUAS.

**Proof.** Consider the non-negative function $L_1 = 0.5(d^2\bar{x}^2 + \bar{v}^2)$, its derivative is calculated as

$$\dot{L}_1 = -\frac{d^3}{c}\bar{x}^2 - d\bar{v}^2 + d^2 D_1 \bar{x} + D_2 \bar{v}$$
$$\leq -c_1 L_1 + c_2(t)\sqrt{L_1}$$

where $c_1 = \frac{2\min\{d^3/c, d\}}{\max\{d^2, 1\}} > 0, c_2(t) = \frac{\max\{d^2|D_1|, |D_2|\}}{\sqrt{0.5\max\{d^2, 1\}}} \geq 0$.

Let $W_1 = \sqrt{L_1} = \sqrt{0.5(d^2\bar{x}^2 + \bar{v}^2)}$, then

$$\dot{W}_1 = \frac{1}{2}\frac{\dot{L}_1}{\sqrt{L_1}} \leq -\frac{1}{2}c_1\sqrt{L_1} + \frac{1}{2}c_2(t) = -\frac{1}{2}c_1 W_1 + \frac{1}{2}c_2(t)$$

Apply the comparison principle [23], we have

$$W_1(t) \leq e^{-\frac{1}{2}c_1 t}W_1(t_0) + \frac{1}{2}\int_{t_0}^{t} e^{-\frac{1}{2}c_1(t-\tau)} c_2(\tau)\,d\tau \qquad (15)$$

In what follows, we denote $\xi = \begin{bmatrix} \xi_1 \\ \xi_2 \end{bmatrix}, \xi_1 = \begin{bmatrix} \bar{x} \\ \bar{v} \end{bmatrix}, \xi_2 = [z, \psi, \bar{u}, r]^T$, and set $t_0 = 0$ for brevity.

As the subsystem (13) is GUAS, so that there exists a class $\mathcal{KL}$ function $\chi_1(\|\xi_2(0)\|_2, t)$ such that $\|\xi_2(t)\|_2 \leq \chi_1(\|\xi_2(0)\|_2, t)$, therefore,

$$c_2(t) = \frac{\max\{d^2|D_1|,|D_2|\}}{\sqrt{0.5\max\{d^2,1\}}} = \frac{\max\left\{d^2\left|\frac{\bar{u}}{c}+\frac{zr}{d}\right|,|-\bar{u}r+d(ar+b\psi)|\right\}}{\sqrt{0.5\max\{d^2,1\}}}$$

$$\leq \frac{1}{\sqrt{0.5\max\{d^2,1\}}}\max\left\{\frac{d^2}{c}+d\chi_1\left(\|\xi_2(0)\|_2,t\right), \chi_1\left(\|\xi_2(0)\|_2,t\right)+d(|a|+|b|)\right\}\chi_1\left(\|\xi_2(0)\|_2,t\right) \quad (16)$$

$$\triangleq \chi_2\left(\|\xi_2(0)\|_2,t\right)$$

Substituting (16) into (15) results

$$W_1(t) \leq e^{-\frac{1}{2}c_1 t}W_1(0) + \frac{1}{2}\chi_2\left(\|\xi_2(0)\|_2,t\right)e^{-\frac{1}{2}c_1 t}\int_0^t e^{\frac{1}{2}c_1\tau}\,d\tau$$

$$= e^{-\frac{1}{2}c_1 t}\sqrt{0.5(d^2\bar{x}^2(0)+\bar{v}^2(0))} + \frac{1}{c_1}\chi_2\left(\|\xi_2(0)\|_2,t\right)\left(1-e^{-\frac{1}{2}c_1 t}\right)$$

$$\leq e^{-\frac{1}{2}c_1 t}\sqrt{0.5(d^2\bar{x}^2(0)+\bar{v}^2(0))} + \frac{1}{c_1}\chi_2\left(\|\xi_2(0)\|_2,t\right)$$

It is hence that

$$\|\xi_1\|_2 \leq \frac{W_1}{\sqrt{0.5\min\{d^2,1\}}} \leq \frac{1}{\sqrt{0.5\min\{d^2,1\}}}\left(\sqrt{0.5(d^2(\bar{x}(0))^2+(\bar{v}(0))^2)}e^{-\frac{1}{2}c_1 t}+\frac{1}{c_1}\chi_2\left(\|\xi_2(0)\|_2,t\right)\right)$$

$$\leq \frac{1}{\sqrt{\min\{d^2,1\}}}\left(\sqrt{0.5\max\{d^2,1\}}\|\xi_1(0)\|_2\, e^{-\frac{1}{2}c_1 t}+\frac{1}{c_1}\chi_2\left(\|\xi_2(0)\|_2,t\right)\right)$$

$$\leq \frac{\sqrt{\max\{d^2,1\}}}{\sqrt{\min\{d^2,1\}}}\|\xi(0)\|_2\, e^{-\frac{1}{2}c_1 t}+\frac{1}{c_1\sqrt{\min\{d^2,1\}}}\chi_2\left(\|\xi(0)\|_2,t\right)$$

$$\triangleq \chi_3\left(\|\xi(0)\|_2,t\right),$$

$$\|\xi(t)\|_2 \leq \|\xi_1(t)\|_2 + \|\xi_2(t)\|_2 \leq \chi_1\left(\|\xi(0)\|_2,t\right) + \chi_3\left(\|\xi(0)\|_2,t\right).$$

Since $\chi_1\left(\|\xi(0)\|_2,t\right), \chi_3\left(\|\xi(0)\|_2,t\right)$ are both class $\mathcal{KL}$ functions of $\left(\|\xi(0)\|_2,t\right)$, so is $\chi_1\left(\|\xi(0)\|,t\right)+\chi_3\left(\|\xi(0)\|,t\right)$, that is, the whole system (12)-(13) is GUAS. This ends the proof of Lemma 1.

□

Now we turn to the controller design of the reduced-order subsystem (13).

**Theorem 1.** System (13) is GUAS under the following control law

$$\begin{aligned}\bar{\tau}_1 &= k_1 zr - k_2\bar{u}, \\ \tau_2 &= -k_3\psi - k_4 r + f(z)\cos(t)\end{aligned} \quad (17)$$

where $k_i\,(i=1,2,3,4)$ are all positive constants, and $f(z)$ a smooth function of $z$ satisfying $f(z)=0 \Leftrightarrow z=0$.

**Proof.** Consider the non-negative function $L_2 = 0.5(k_1 z^2 + \bar{u}^2)$, its derivative along the solution of the closed-loop system (13) (17) is

$$\dot{L}_2 = -k_1 z \bar{u} r + \bar{u} \tau_1 = -k_2 \bar{u}^2 \leq 0,$$

which implies that, $(L_2, z, \bar{u})$ are all bounded and $L_2$ converges to a constant limit $L_2(\infty)$ as $t \to \infty$.

Since $z$ is bounded, so that $(\psi, r, \dot{r})$ are all well defined and bounded in view of the closed-loop dynamics $\dot{\psi} = r, \dot{r} = -k_3 \psi - k_4 r + f(z)\cos(t)$.

As $\dot{\bar{u}} = \tau_1 = k_1 z r - k_2 \bar{u} \in L_\infty$, so is $\ddot{L}_2$, we thus conclude $\lim_{t \to \infty} \bar{u} = 0$ by Barbalat Lemma [23].

As $(\bar{u}, r, z, \tau_1, \tau_2 = \dot{r}) \in L_\infty$, so is $\ddot{\bar{u}} = \dot{\tau}_1 = k_1(-\bar{u}r)r + k_1 z \tau_2 - k_2 \tau_1 \in L_\infty$, and hence $\lim_{t \to \infty} \dot{\bar{u}} = \tau_1 = 0$ can be induced, which in turn implies $\lim_{t \to \infty}(zr) = 0$ from the expression of $\dot{\bar{u}}$.

By $\lim_{t \to \infty} \bar{u} = 0$ and $\lim_{t \to \infty} L_2 = 0.5 \lim_{t \to \infty}(k_1 z^2 + \bar{u}^2) = L_2(\infty)$, we know that $z$ converges to a constant limit $z(\infty)$, therefore, we can conclude $z(\infty) \lim_{t \to \infty} r = 0$ from $\lim_{t \to \infty}(zr) = 0$.

If $z(\infty) \neq 0$, then $\lim_{t \to \infty} r = 0$. It can be easily verified that $\ddot{r} \in L_\infty$, $r^{(3)} \in L_\infty$, we thus have $\lim_{t \to \infty} \dot{r} = \lim_{t \to \infty} \ddot{r} = 0$ by Barbalat Lemma [23]. By the expression of 

$$\ddot{r} = \dot{\tau}_2 = -k_3 r - k_4 \dot{r} + \frac{\partial f(z)}{\partial z}(-r\bar{u})\cos(t) - f(z)\sin(t),$$ 

one gets $\lim_{t \to \infty}(f(z)\sin(t)) = 0$, implying $\lim_{t \to \infty}(f(z(\infty))\sin(t)) = 0$, this is, $f(z(\infty)) = 0$ and thus $z(\infty) = 0$, which clearly contradicts the assumption of $z(\infty) \neq 0$, therefore, $z(\infty)$ must be zero.

By $z(\infty) = 0$, we conclude $\lim_{t \to \infty}(\psi, r) = 0$ referring to the closed-loop dynamics of $\dot{\psi} = r, \dot{r} = -k_3 \psi - k_4 r + f(z)\cos(t)$.

Since all the states of the closed-loop system (13) (15) are globally uniformly convergent to zero, and the origin is an equilibrium point, so that the closed-loop system (13) (15) is GUAS. This ends the proof of Theorem 1.

□

## 4. Trajectory tracking

In this section, we propose a tracking controller that makes the ship asymptotically track the reference one with an exponential convergence rate via mild PE conditions.

The tracking errors in respect to the ship coordinate frame are defined as

$$\begin{cases} x_e = (x - x_d)\cos\psi + (y - y_d)\sin\psi, \\ y_e = -(x - x_d)\sin\psi + (y - y_d)\cos\psi, \\ \psi_e = \psi - \psi_d, \\ u_e = u - u_d, r_e = r - r_d, v_e = v - v_d, \\ \tau_{1e} = \tau_1 - \tau_{1d}, \tau_{2e} = \tau_2 - \tau_{2d} \end{cases} \quad (18)$$

The dynamics of the tracking error can then be derived from (3), (4), (7), (8), (18) as

$$\begin{aligned}
\dot{x}_e &= u_e - u_d(\cos\psi_e - 1) - v_d \sin\psi_e + r y_e \\
&= u_e - u_d(\cos\psi_e - 1 + 0.5\psi_e^2) - v_d(\sin\psi_e - \psi_e) + 0.5 u_d \psi_e^2 - v_d \psi_e + r y_e, \\
&= u_e - \alpha\psi_e + 0.5 u_d \psi_e^2 - v_d \psi_e + r y_e \\
\dot{y}_e &= v_e - v_d(\cos\psi_e - 1) + u_d \sin\psi_e - r x_e \\
&= v_e - v_d(\cos\psi_e - 1 + 0.5\psi_e^2) + u_d(\sin\psi_e - \psi_e) + 0.5 v_d \psi_e^2 + u_d \psi_e - r x_e \\
&= v_e + \beta\psi_e + 0.5 v_d \psi_e^2 + u_d \psi_e - r x_e, \\
\dot{\psi}_e &= r_e, \\
\dot{v}_e &= -a\tau_{2e} - b r_e - c u_e r - c u_d r_e - d v_e, \\
\dot{r}_e &= \tau_{2e}, \\
\dot{u}_e &= \tau_{1e}
\end{aligned} \quad (19)$$

where

$$\alpha \triangleq \begin{cases} \alpha^* = \dfrac{u_d(\cos\psi_e - 1 + 0.5\psi_e^2) + v_d(\sin\psi_e - \psi_e)}{\psi_e} & \psi_e \neq 0 \\ 0 & \psi_e = 0 \end{cases}$$

$$\beta \triangleq \begin{cases} \beta^* = \dfrac{-v_d(\cos\psi_e - 1 + 0.5\psi_e^2) + u_d(\sin\psi_e - \psi_e)}{\psi_e} & \psi_e \neq 0 \\ 0 & \psi_e = 0 \end{cases}$$

$$\dot{\alpha} \triangleq \begin{cases} \dfrac{d\alpha^*}{dt} & \psi_e \neq 0 \\ 0 & \psi_e = 0 \end{cases}, \dot{\beta} \triangleq \begin{cases} \dfrac{d\beta^*}{dt} & \psi_e \neq 0 \\ 0 & \psi_e = 0 \end{cases}$$

It can be easily checked that $\lim\limits_{\psi_e \to 0} \alpha^* = \lim\limits_{\psi_e \to 0} \dot{\alpha}^* = \lim\limits_{\psi_e \to 0} \beta^* = \lim\limits_{\psi_e \to 0} \dot{\beta}^* = 0$, thus $\alpha, \dot{\alpha}, \beta, \dot{\beta}$ are all continuous functions of $\psi_e$.

In order to remove the terms $-a\tau_{2e} - b r_e - c u_d r_e$ in the $v_e$ dynamics (forth equation of (19)) and the term $v_e$ in the $y_e$ dynamics (second equation of (19)), we introduce the new error state and input variables as follows

$$\begin{aligned}
\bar{v}_e &= v_e + ar_e + b\psi_e + cu_d\psi_e, \\
z_e &= dy_e + \bar{v}_e, \\
\bar{u}_e &= dx_e + cu_e, \\
\bar{\tau}_{1e} &= \dot{\bar{u}}_e = d(u_e - u_d(\cos\psi_e - 1) - v_d\sin\psi_e + ry_e) + c\tau_{1e}
\end{aligned} \quad (20)$$

The new error dynamics can be deduced from (19) and (20) as

$$\begin{aligned}
\dot{\bar{v}}_e &= -a\tau_{2e} - br_e - (\bar{u}_e - dx_e)r - cu_d r_e - dv_e + a\tau_{2e} + br_e + c\dot{u}_d\psi_e + cu_d r_e \\
&= -dv_e - (\bar{u}_e - dx_e)r + c\dot{u}_d\psi_e \\
&= -d(\bar{v}_e - ar_e - b\psi_e - cu_d\psi_e) - (\bar{u}_e - dx_e)r + c\dot{u}_d\psi_e, \\
\dot{z}_e &= dv_e + d\beta\psi_e + d(0.5v_d\psi_e^2 + u_d\psi_e) - drx_e - dv_e - (\bar{u}_e - dx_e)r + c\dot{u}_d\psi_e \\
&= (d\beta + c\dot{u}_d + 0.5dv_d\psi_e + du_d)\psi_e - \bar{u}_e r, \\
\dot{x}_e &= \frac{\bar{u}_e - dx_e}{c} - \alpha\psi_e + 0.5u_d\psi_e^2 - v_d\psi_e + r\frac{z_e - \bar{v}_e}{d}, \\
\dot{\psi}_e &= r_e, \\
\dot{r}_e &= \tau_{2e}, \\
\dot{\bar{u}}_e &= \bar{\tau}_{1e}
\end{aligned}$$

The above error state dynamics can be divided into the following two subsystems

$$\begin{aligned}
\dot{x}_e &= -\frac{d}{c}x_e - \frac{1}{d}(r_e + r_d)\bar{v}_e + D_3, \\
\dot{\bar{v}}_e &= -d\bar{v}_e + dx_e(r_e + r_d) + D_4,
\end{aligned} \quad (21)$$

$$\begin{aligned}
\dot{z}_e &= (d\beta + c\dot{u}_d + 0.5dv_d\psi_e + du_d)\psi_e - \bar{u}_e(r_e + r_d), \\
\dot{\psi}_e &= r_e, \\
\dot{r}_e &= \tau_{2e}, \\
\dot{\bar{u}}_e &= \bar{\tau}_{1e}
\end{aligned} \quad (22)$$

where

$$\begin{aligned}
D_3 &= \frac{\bar{u}_e}{c} - \alpha\psi_e + 0.5u_d\psi_e^2 - v_d\psi_e + \frac{z_e(r_e + r_d)}{d}, \\
D_4 &= d(ar_e + b\psi_e + cu_d\psi_e) - \bar{u}_e(r_e + r_d) + c\dot{u}_d\psi_e
\end{aligned} \quad (23)$$

**Remark 2.** Observe that system (21)-(22) enjoys a cascaded structure such that the stabilization problem of the whole system (21)-(22) can be reduced to the one of its subsystem (22) as claimed in the following Lemma.

**Lemma 2.** Suppose that $(u_d, r_d, \dot{u}_d, \dot{r}_d)$ are uniformly bounded, then any control law $\bar{\tau}_{1e}(t, z_e, \psi_e, \bar{u}_e, r_e), \tau_{2e}(t, z_e, \psi_e, \bar{u}_e, r_e)$ that makes the closed-loop subsystem (22) GKES also makes the whole closed-loop system (21)-(22) GKES.

**Proof.** Consider the non-negative function $L_2 = 0.5(d^2 x_e^2 + \bar{v}_e^2)$, its derivative is calculated as

$$\dot{L}_2 = -\frac{d^3}{c} x_e^2 - d\bar{v}_e^2 + d^2 D_3 x_e + D_4 \bar{v}_e \leq -c_3 L_2 + c_4(t)\sqrt{L_2}$$

where $c_3 = \frac{2\min\{d^3/c, d\}}{\max\{d^2, 1\}} > 0, c_4(t) = \frac{\max\{d^2 |D_3|, |D_4|\}}{\sqrt{0.5\max\{d^2, 1\}}} \geq 0$.

Let $W_2 = \sqrt{L_2} = \sqrt{0.5(d^2 x_e^2 + \bar{v}_e^2)}$, then

$$\dot{W}_2 = \frac{1}{2} \frac{\dot{L}_2}{\sqrt{L_2}} \leq -\frac{1}{2} c_3 \sqrt{L_2} + \frac{1}{2} c_4(t) = -\frac{1}{2} c_3 W_2 + \frac{1}{2} c_4(t)$$

Apply the comparison principle [23], we have

$$W_2(t) \leq e^{-\frac{1}{2}c_3 t} W_2(t_0) + \frac{1}{2} \int_{t_0}^{t} e^{-\frac{1}{2}c_3(t-\tau)} c_4(\tau) d\tau \qquad (24)$$

As $(u_d, r_d, \dot{r}_d)$ are uniformly bounded, we conclude $(v_d, \dot{v}_d)$ are also uniformly bounded in view of $\dot{v}_d = -a\dot{r}_d - br_d - cu_d r_d - dv_d$.

In what follows, we denote $\eta = \begin{bmatrix} \eta_1 \\ \eta_2 \end{bmatrix}$, $\eta_1 = \begin{bmatrix} x_e \\ \bar{v}_e \end{bmatrix}$, $\eta_2 = [z_e, \psi_e, \bar{u}_e, r_e]^T$, and set $t_0 = 0$ for brevity.

Since the closed-loop system (22) is GKES, then there exists a class $\kappa$ function $\rho_1(\cdot)$ and a sufficient small positive constant $0 < \gamma < 2c_3$ such that $\|\eta_2(t)\|_2 \leq \rho_1\left(\|\eta_2(0)\|_2\right) e^{-\gamma t}$.

In view of the assumption $(u_d, \dot{u}_d, r_d, v_d) \in L_\infty$ and the expression of $c_4(t)$, it can be easily verified that there exists another class $\kappa$ function $\rho_2(\cdot)$ such that

$$c_4(t) \leq \rho_2\left(\|\eta_2(0)\|_2\right) e^{-\gamma t} \qquad (25)$$

Substituting (25) into (24) results

$$W_2(t) \le e^{-\frac{1}{2}c_3 t} W_2(0) + \frac{1}{2}\rho_2\left(\|\eta_2(0)\|_2\right) e^{-\frac{1}{2}c_3 t} \int_0^t e^{(\frac{1}{2}c_3 - \gamma)\tau} \, d\tau$$

$$= e^{-\frac{1}{2}c_3 t} \sqrt{0.5(d^2 x_e^2(0) + \bar{v}_e^2(0))} + \rho_2\left(\|\eta_2(0)\|_2\right) \frac{1}{c_3 - 2\gamma}\left(e^{-\gamma t} - e^{-\frac{1}{2}c_3 t}\right)$$

$$\le e^{-\frac{1}{2}c_3 t} \sqrt{0.5(d^2 x_e^2(0) + \bar{v}_e^2(0))} + \frac{1}{c_3 - 2\gamma} \rho_2\left(\|\eta_2(0)\|_2\right) e^{-\gamma t}$$

$$\le \left(\sqrt{0.5\max\{d^2,1\}}\|\eta_1(0)\|_2 + \frac{1}{c_3 - 2\gamma}\rho_2\left(\|\eta_2(0)\|_2\right)\right) e^{-\gamma t}$$

$$\le \left(\sqrt{0.5\max\{d^2,1\}}\|\eta(0)\|_2 + \frac{1}{c_3 - 2\gamma}\rho_2\left(\|\eta(0)\|_2\right)\right) e^{-\gamma t}$$

$$\triangleq \rho_3\left(\|\eta(0)\|_2\right) e^{-\gamma t}$$

It is hence that

$$\|\eta_1\|_2 \le \frac{W_2}{\sqrt{0.5\min\{d^2,1\}}} \le \frac{\rho_3\left(\|\eta(0)\|_2\right) e^{-\gamma t}}{\sqrt{0.5\min\{d^2,1\}}},$$

$$\|\eta\|_2 \le \|\eta_1\|_2 + \|\eta_2\|_2 \le \frac{\rho_3\left(\|\eta(0)\|_2\right) e^{-\gamma t}}{\sqrt{0.5\min\{d^2,1\}}} + \rho_1\left(\|\eta(0)\|_2\right) e^{-\gamma t}$$

$$= \left(\frac{\rho_3\left(\|\eta(0)\|_2\right)}{\sqrt{0.5\min\{d^2,1\}}} + \rho_1\left(\|\eta(0)\|_2\right)\right) e^{-\gamma t} \triangleq \rho_4\left(\|\eta(0)\|_2\right) e^{-\gamma t}$$

As $\rho_1(\cdot), \rho_3(\cdot)$ are both class $\kappa$ functions, so is $\rho_4(\cdot)$. Therefore, the closed-loop system (21)-(22) is GKES. This completes the proof of Lemma 2.

□

Now we turn to the controller design for subsystem (22).

**Theorem 2.** Suppose that the reference surge and yaw velocities $(u_d, r_d)$ satisfy the following persistent exciting conditions

$$(u_d, \dot{u}_d, \ddot{u}_d, r_d, \dot{r}_d) \in L_\infty, \lim_{t \to \infty}\left(|u_d| + |r_d|\right) \ne 0 \tag{26}$$

then the following control law

$$\begin{aligned} r_{ed} &= -k_1 z_e (d\beta + c\dot{u}_d + 0.5 dv_d \psi_e + du_d) - k_2 \psi_e, \\ \bar{\tau}_{1e} &= k_1 z_e r - k_3 \bar{u}_e, \\ \tau_{2e} &= \dot{r}_{ed} - \psi_e - k_4(r_e - r_{ed}) \end{aligned} \tag{27}$$

makes the origin of the closed-loop system (22) (27) GKES, where $(k_1, k_2, k_3, k_4)$ are all positive constants, and $r_{ed}$ is an assistant variable.

**Proof.** Consider the following positive function

$$L_3 = 0.5(k_1 z_e^2 + \psi_e^2 + (r_e - r_{ed})^2 + \overline{u}_e^2)$$

Its derivative along the closed-loop system (22) (27) is calculated as

$$\begin{aligned}
\dot{L}_3 &= k_1 z_e (d\beta + c\dot{u}_d + 0.5 dv_d \psi_e + du_d)\psi_e - k_1 z_e \overline{u}_e r + \psi_e r_e + (r_e - r_{ed})(\tau_{2e} - \dot{r}_{ed}) + \overline{u}_e \overline{\tau}_{1e} \\
&= k_1 z_e (d\beta + c\dot{u}_d + 0.5 dv_d \psi_e + du_d)\psi_e - k_1 z_e \overline{u}_e r + \psi_e (r_e - r_{ed}) + r_{ed} \psi_e \\
&\quad + (r_e - r_{ed})(\tau_{2e} - \dot{r}_{ed}) + \overline{u}_e \overline{\tau}_{1e} \\
&= \left[ k_1 z_e (d\beta + c\dot{u}_d + 0.5 dv_d \psi_e + du_d) + r_{ed} \right] \psi_e + \overline{u}_e (-k_1 z_e r + \overline{\tau}_{1e}) + (r_e - r_{ed})(\tau_{2e} - \dot{r}_{ed} + \psi_e) \\
&= -k_2 \psi_e^2 - k_3 \overline{u}_e^2 - k_4 (r_e - r_{ed})^2 \leq 0
\end{aligned}$$

Hence $L_3$ is uniformly bounded, non-increasing, convergent to a constant limit $L_3(\infty)$ as $t \to \infty$, and the closed-loop system (22)-(27) is globally uniformly stable.

From $\dot{v}_d = -a\dot{r}_d - br_d - cu_d r_d - dv_d$, and $(u_d, r_d, \dot{r}_d) \in L_\infty$, we have $(v_d, \dot{v}_d) \in L_\infty$. In view of $(u_d, \dot{u}_d, \ddot{u}_d, r_d, \dot{r}_d, v_d, \dot{v}_d) \in L_\infty$, we can deduce that all the state variables and their derivatives are bounded.

By verifying $\ddot{L}_3 \in L_\infty$, one obtains $\lim_{t \to \infty} \dot{L}_3 = 0 \Rightarrow \lim_{t \to \infty}(\psi_e, \overline{u}_e, r_e - r_{ed}) = 0$.

As $\ddot{\psi}_e = \tau_{2e} \in L_\infty$, we conclude $\lim_{t \to \infty} r_e = \lim_{t \to \infty} \dot{\psi}_e = 0 \Rightarrow \lim_{t \to \infty} r_{ed} = \lim_{t \to \infty}(r_e - (r_e - r_{ed})) = 0$,

so that $\lim_{t \to \infty}\left(-k_1 z_e (d\beta + c\dot{u}_d + 0.5 dv_d \psi_e + du_d) - k_2 \psi_e\right) = 0 \Rightarrow \lim_{t \to \infty}(z_e (c\dot{u}_d + du_d)) = 0$.

Since $(L_3, \psi_e, \overline{u}_e, r_e - r_{ed})$ have constant limits as $t \to \infty$, so that $z_e$ converges to a constant limit $z_e(\infty)$, implying that, $z_e(\infty) \lim_{t \to \infty}(du_d + c\dot{u}_d) = 0$.

On the other hand, we can also infer $z_e(\infty) \lim_{t \to \infty} r_d = 0$ from $\lim_{t \to \infty} \overline{u}_e = 0, \ddot{\overline{u}}_e = \dot{\overline{\tau}}_{1e} \in L_\infty \Rightarrow$

$\lim_{t \to \infty} \overline{\tau}_{1e} = \lim_{t \to \infty}(k_1 z_e r - k_3 \overline{u}_e) = 0 \Rightarrow \lim_{t \to \infty}(z_e r) = 0 \Rightarrow z_e(\infty) \lim_{t \to \infty} r_d = 0$.

If $z(\infty) \neq 0$, we have $\lim_{t \to \infty}(du_d + c\dot{u}_d) = 0 \Rightarrow \lim_{t \to \infty} u_d = 0$ and $\lim_{t \to \infty} r_d = 0$, which clearly contradicts the persistent exciting condition (26), it is hence that $z(\infty)$ must be zero.

Up to now, we have shown that all the states of the closed-loop system (22) (27) are globally uniformly convergent to zero, and the origin is globally uniformly stable, we finally conclude that the closed-loop system (22) (27) is GUAS.

The linearization of the closed-loop system (22) (27) at the origin can be written as

$$\begin{aligned}
\dot{z}_e &= (c\dot{u}_d + du_d)\psi_e - \bar{u}_e r_d, \\
\dot{\psi}_e &= r_e, \\
\dot{r}_e &= \tau_{2e}, \\
\dot{\bar{u}}_e &= \bar{\tau}_{1e}, \\
r_{ed} &= -k_1 z_e (c\dot{u}_d + du_d) - k_2 \psi_e, \\
\bar{\tau}_{1e} &= k_1 z_e r - k_3 \bar{u}_e, \\
\tau_{2e} &= \dot{r}_{ed} - \psi_e - k_4 (r_e - r_{ed})
\end{aligned} \qquad (28)$$

Take the positive function $L_4 = 0.5(k_1 z_e^2 + \psi_e^2 + (r_e - r_{ed})^2 + \bar{u}_e^2)$, its derivative along (28) is $\dot{L}_4 = -k_2 \psi_e^2 - k_3 \bar{u}_e^2 - k_4 (r_e - r_{ed})^2 \leq 0$. Along the same line of the proof for closed-loop system (22) (27), one can show that system (28) is also GUAS.

Since the closed-loop nonlinear system (22) (27) and its linearization (28) are both GUAS, and the Jacobian matrix of (22) is uniformly bounded, it is thus GKES [25]. This ends the proof of Theorem 2.

□

**Remark 3.** The PE condition (26) is simple, easy for checking, and weak as what required is that the reference ship keeps moving instead of standing still or converging to a fixed point.

## 5. Simulations

In this section, the effectiveness of the proposed control laws is verified via simulation examples. Consider an underactuated surface ship with model parameters as [19]

$$M = \begin{bmatrix} 25.8 & 0 & 0 \\ 0 & 33.8 & 1.0115 \\ 0 & 1.0115 & 2.76 \end{bmatrix}, D = \begin{bmatrix} 0.9257 & 0 & 0 \\ 0 & 2.8909 & -0.2601 \\ 0 & -0.2601 & 0.5 \end{bmatrix}.$$

For point stabilization, the control parameters are assigned to $k_1 = 0.6, k_2 = 0.4, k_3 = 0.1, k_4 = 0.1$; $f(z)$ is selected as $f(z) = 10\tanh(10z^2)$. Simulation results for the initial state $(x, y, \psi, u, v, r)_0 = (-2, 2, 0, 0, 0, 0)$ and $(x, y, \psi, u, v, r)_0 = (0, 2, 0, 0, 0, 0)$ are shown in Fig.1-Fig-2 respectively.

For trajectory tracking, the reference trajectories are taken as a straight-line one generated by $\tau_{1d} = \tau_{2d} = 0$, $(x_d, y_d, \psi_d, u_d, v_d, r_d)_0 = (0, 0, \pi/8, 4, 0, 0)$ and a circular one generated by $\tau_{1d} = \tau_{2d} = 0, (x_d, y_d, \psi_d, u_d, v_d, r_d)_0 = (-2, 1, 0.2, -0.32, 0.188)$. The control parameters are assigned to $k_1 = 1, k_2 = 0.5, k_3 = 0.5, k_4 = 1$. Simulation results are shown in Fig.3-4 respectively.

It is observed from Fig.1-Fig.4 that the proposed control laws successfully steer the state to the desired state and the desired reference trajectories.

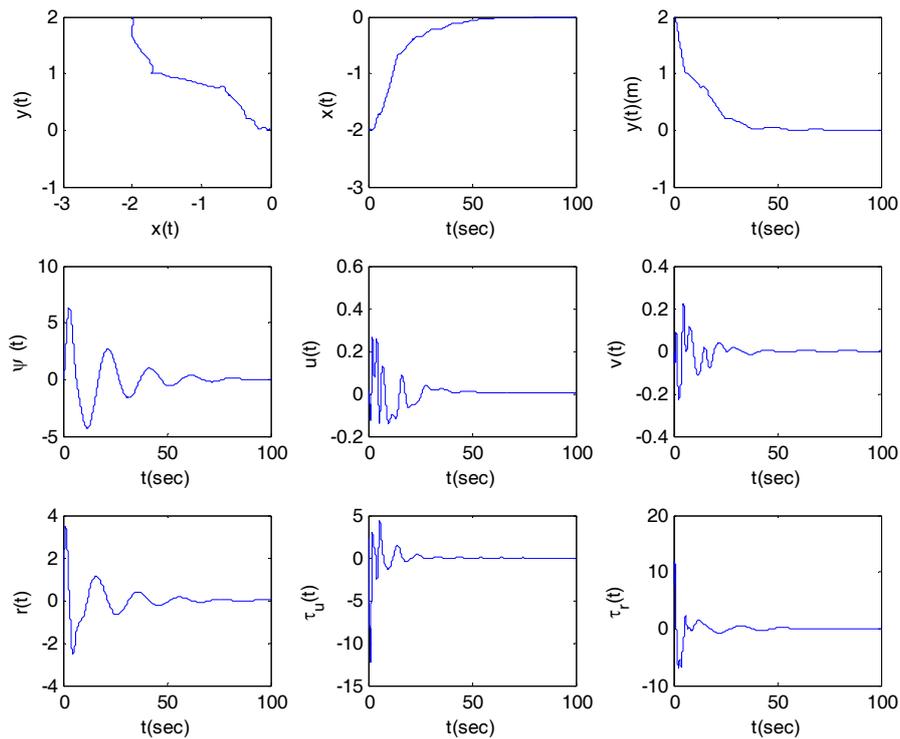

**Fig.1** Plots of the geometric path, the time trajectories of states and control inputs for point stabilization with an initial state $(x, y, \psi, u, v, r)_0 = (-2, 2, 0, 0, 0, 0)$.

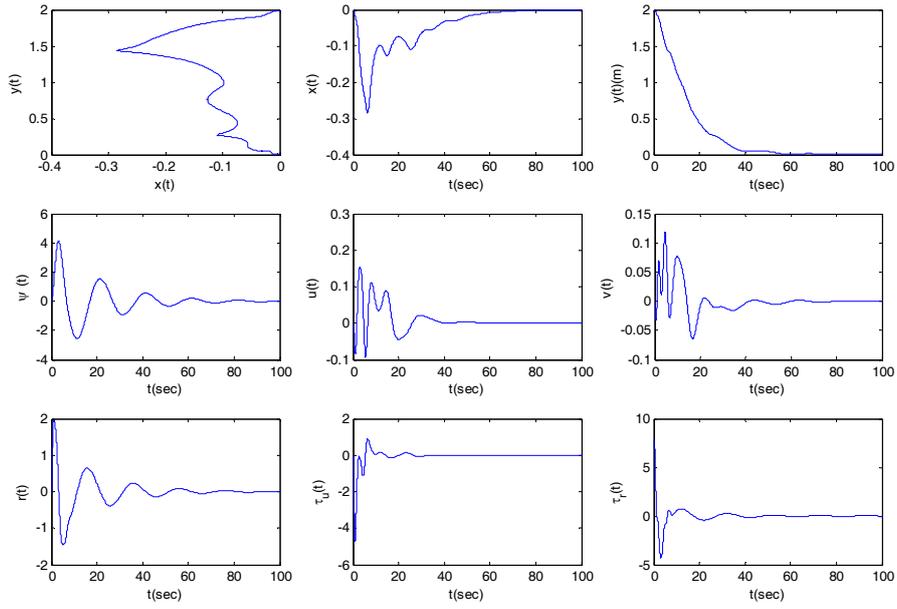

**Fig.2** Plots of the geometric path, the time trajectories of states and control inputs for point stabilization with an initial state $(x, y, \psi, u, v, r)_0 = (0, 2, 0, 0, 0, 0)$.

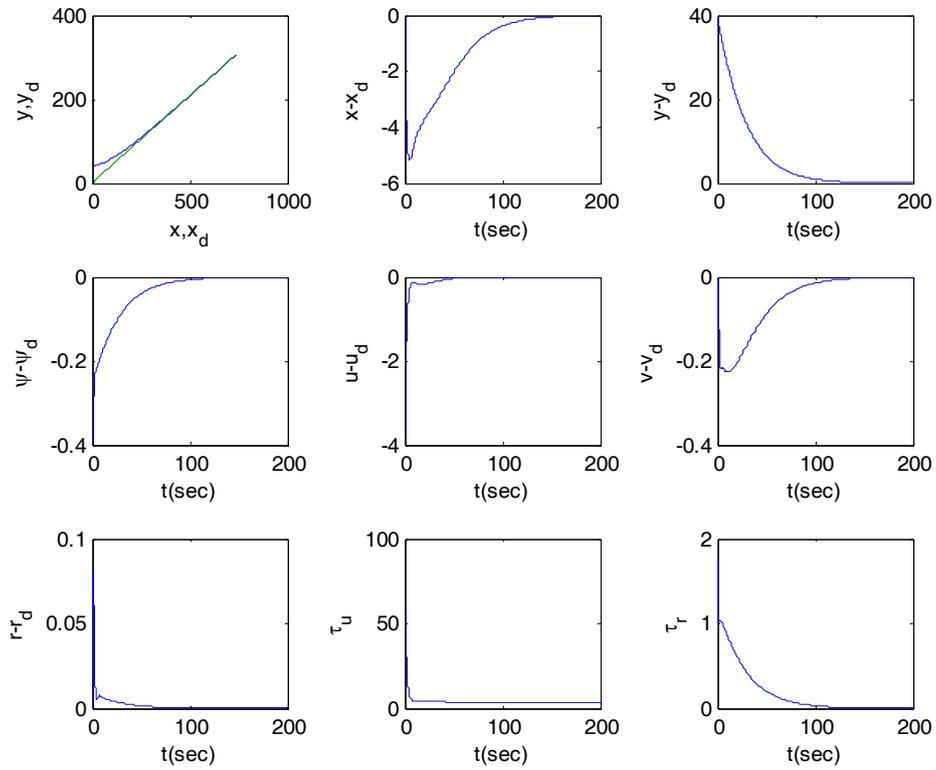

**Fig.3** Plots of the reference and actual geometric paths, the time history of state tracking errors

and control inputs for straight-line trajectory tracking with an initial state $(x, y, \psi, u, v, r)_0 = (0, 40, 0, 0, 0, 0)$.

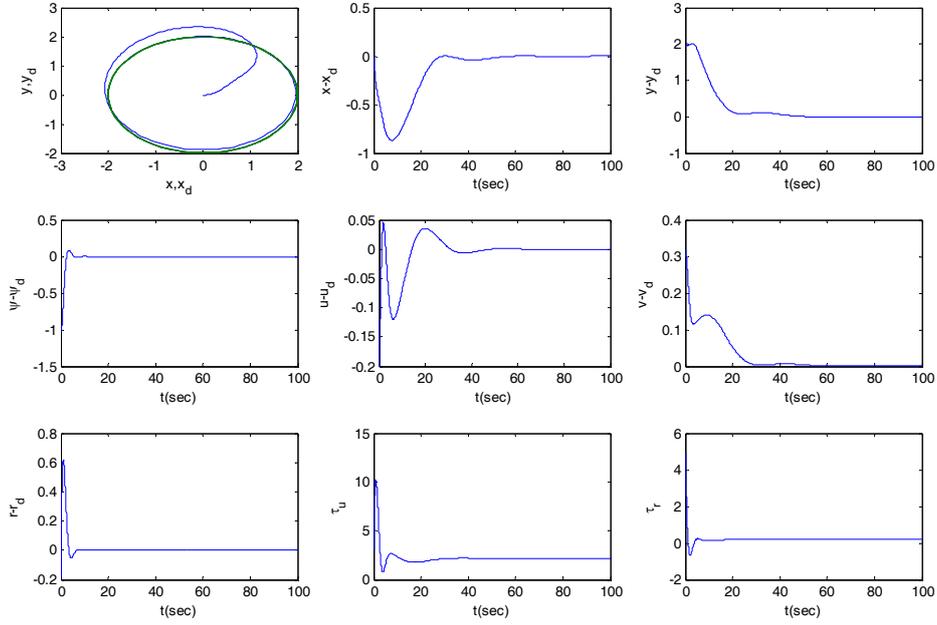

**Fig.4** Plots of the reference and actual geometric paths, the time history of state tracking errors and control inputs for circular trajectory tracking with the zero initial state.

## 6. Conclusion

In this work, we have developed full state stabilization and full state trajectory tracking control schemes for underactuated surface ships with non-diagonal inertia and damping matrices. The proposed control laws are the first ones guaranteeing global uniform asymptotic convergence of state to the desired set point and global $\kappa-$exponential convergence to the desired reference trajectory for such a ship model. Furthermore, the proposed trajectory tracking controller requires simple and mild persistent exciting conditions.

It is noted that the construction of the control laws rely on the exact ship model with no parameter uncertainties and no external disturbances. Our future research topics will focus on designing robust stabilizing and tracking controllers for underactuated surface ships with unknown non-diagonal inertia/damping matrices and unknown external disturbances.

**Acknowledgments**

This work was supported by National Science Foundation of China under grants 60874012.

**References**

[1] T.I.Fossen, Marine Control Systems, Trondheim Norway: Marine Cybernetics, 2002.

[2] K.Y.Pettersen and O.Egeland, Exponential stabilization of an underactuated surface vessels, Proc.35$^{th}$ IEEE Conference on Decision and Control, Kobe, Japan, December 1996, pp.967-972.

[3] M.Reyhanoglu, Exponential stabilization of an underactuated autonomous surface vessel, Automatica, 1997, 33(12): 2249-2254.

[4] F.Mazenc, K.Y.Pettersen, H.Nijmeijer, Global uniform asymptotic stabilization of an underactuated surface vessel, IEEE Transactions on Automatic Control, 2002, 47(10):1759-1762.

[5] Wenjie Dong, Yi Guo, Global time-varying stabilization of underactuated surface vessel, IEEE Transactions on Automatic Control, 2005, 50(6):859-864.

[6] Ma Bao-Li, Global k-exponential asymptotic stabilization of underactuated surface vessels, Systems and Control Letters, 2009, 58(3): 194-201.

[7] Ma Bao-Li, Huo Wei, Smooth time-varying uniform asymptotic stabilization of underactuated surface vessels, Joint 48$^{th}$ IEEE Conference on Decision and Control, Shanghai, P.R.China, Dec.16-18, 2009, pp.3137-3141.

[8] J.Ghommam, F.Mnif, N.Derbel, Global stabilization and tracking control of underactuated surface vessels, IET Control Theory and Applications, 2010, 4(1): 71-78.

[9] K.Y.Pettersen and H.Nijmeijer, Underactuated ship tracking control: theory and experiments, International Journal of Control, 2001, 74(14):1435-1446.

[10] K.D.Do, Z.P.Jiang, and J.Pan, Underactuated ship global tracking under relaxed conditions, IEEE Transactions on Automatic Control, 2002, 47(9):1529-1536.

[11] E.Lefeber, K.Y.Pettersen, and H.Nijmeijer, Tracking control of an underactuated ship, IEEE Transactions on Control Systems Technology, 2003, 11(1): 52-61.

[12] Ti-Chung Lee and Zhong-Ping Jiang, New cascade approach for global k-exponential tracking of underactuated ships, IEEE Transactions on Automatic Control, 2004, .49(12): 2297-2203.

[13] K.D.Do, Z.P.Jiang, J.Pan, Universal controllers for stabilization and tracking of underactuated ships, Systems and Control Letters, 2002, 47(4): 299-317.

[14] A.Behal, D.M.Dawson, W.E.Dixon, and Y.Fang, Tracking and regulation control of an underactuated surface vessel with non-integrable dynamics, IEEE Transactions on Automatic Control, 2002, 47(3):495-500.

[15] A.P.Aguiar, A.M.Pascoal, Dynamic positioning and way-point tracking of underactuated AUVs in the presence of ocean currents, International Journal of Control, 2007, 80(7): 1092-1108.

[16] K.D.Do, Z.P.Jiang and J.Pan, Robust global stabilization of underactuated ships on a linear course: state and output feedback, International Journal of Control, 2003, 76(1):1-17.

[17] K.D.Do, J.Pan, State- and output-feedback robust path-following controllers for underactuated ships using Serret-Frenet frame, Ocean Engineering, 2004, 31(16): 1967-1997.

[18] K.D.Do, J.Pan, Global robust adaptive path following of underactuated ships, Automatica , 2006, 42(10):1713-1722.

[19] E.Fredriksen and K.Y.Pettersen, Global $\kappa-$ exponential way-point maneuvering of ships: Theory and experiments, Automatica, 2006, 42:677-687.

[20] K.D.Do, J.Pan, Global tracking control of underactuated ships with nonzero off-diagonal terms in their system matrices, Automatica, 2005, 41: 87-95.

[21] K.D.Do, J.Pan, Underactuated ships follow smooth paths with integral actions and without velocity measurements for feedback: theory and experiments. IEEE Transactions on Control Systems technology, 2006, 14(2):308-322.

[22] Ji-Hong Li, Pan-Mook Lee, Bong-Huan Jun,Yong-Kon Lim, Point-to-point navigation of underactuated ships, Automatica, 2008, 44: 3201-3205.

[23] H.K.Khalil, Nonlinear Systems (Second-Edition), Prentice-Hall, Upper Saddle River, New Jersey, USA, 1996.

[24] O.J.Sordalen, O.Egeland, Exponential stabilization of nonholonomic chained systems, IEEE Translation on Automatic Control 40 (1995) 35-39.

[25] A.A.J. Lefeber, Tracking control of nonlinear mechanical systems, Ph.D. Dissertation, Dept. Mech. Eng., University of Twente, , Twente, The Netherlands, 2000.